\documentclass[11pt]{article}

\usepackage{amsfonts}
\usepackage{amsmath}
\usepackage{amstext}
\usepackage{amsthm}
\usepackage{xspace}

\input psfig.sty

\newcommand{\N}{\mathbb N}

\newcommand{\E}{\mathbb E}
\newcommand{\A}{{\mathcal A}}
\newcommand{\B}{{\mathcal B}}
\newcommand{\oq}{\overline q}
\newcommand{\uq}{\underline q}

\newtheorem{teorema}{Theorem}[section]
\newtheorem{prop}[teorema]{Proposition}
\newtheorem{lemma}[teorema]{Lemma}
\theoremstyle{definition}
\newtheorem{definizione}[teorema]{Definition}
\newtheorem{cor}[teorema]{Corollary}

\begin{document}

\title{Algorithmic information for intermittent systems with an indifferent fixed point}
\author{Claudio Bonanno\footnote{Dipartimento di Matematica e Informatica,
Universit\`a di Camerino, via Madonna delle Carceri 9, 62032
Camerino (MC), Italy, email $<$claudio.bonanno@unicam.it$>$} \and
Stefano Galatolo\footnote{Dipartimento di Matematica Applicata,
Universit\`a di Pisa, via Bonanno 26/b, 56125 Pisa, Italy, email
$<$galatolo@mail.dm.unipi.it$>$}}

\maketitle

\begin{abstract}
Measuring the average information that is necessary to describe
the behaviour of a dynamical system leads to a generalization of
the Kolmogorov-Sinai entropy. This is particularly interesting
when the system has null entropy and the information increases
less than linearly with respect to time. We consider two classes
of maps of the interval with an indifferent fixed point at the
origin and an infinite natural invariant measure. We calculate
that the average information that is necessary to describe the
behaviour of its orbits increases with time $n$ approximately as
$n^\alpha $, where $\alpha<1$ depends only on the asymptotic
behaviour of the map near the origin.
\end{abstract}

\section{Introduction} \label{sec:intro}

The complexity and unpredictability of a chaotic system has been
measured using many different indicators. Among all one of the
most important is the Kolmogorov-Sinai (K-S) entropy. Being based
on the Shannon's notion of information, it is an average measure
of the quantity of information that is necessary to describe each
step of the behaviour of the system (with an arbitrary accuracy
given by the choice of a partition).

More recently other notions of information content, such as the
Kolmo\-go\-rov-Chaitin \emph{Algorithmic Information Content},
have been applied to dynamical systems. These notions are
pointwise and allow to consider the complexity of the behaviour of
a single orbit. Hence it is possible to define the information
$I(x,n,Z)$ contained in $n$ steps of the orbit of a point $x$ with
respect to a partition $Z$ of the phase space. This can be done by
associating to the orbit of $x$ the symbolic orbit with respect to
$Z$ and considering the information content of this string (see
Section \ref{sec:info}). The average of this pointwise information
over an invariant measure $\mu$ is strictly related to the K-S
entropy $h_\mu(T,Z)$ of the measure $\mu$ relative to the
partition $Z$. Indeed for a "typical" point $x$ it holds
$I(x,n,Z)\sim n\ h_\mu(T,Z)$ (see Theorem \ref{teo:ks-aic}). When
the entropy is null the previous relation becomes $I(x,n,Z)=o(n)$.
The many possible different sublinear asymptotic behaviours of
$I(x,n,Z)$ correspond to different kinds of "weakly" chaotic
dynamics. The importance of this indicator of weak chaos is also
confirmed by the relations that have been proved, even in the null
entropy case, between the behaviour of the information and many
important features of the dynamics, such as sensitivity to initial
conditions, dimensions, recurrence
(\cite{brudno},\cite{white},\cite{bgi},\cite{gal1},\cite{keller})
and global topological complexity indicators (\cite{gal2}, see
also \cite{afr} for relations between the topological complexity
and other physically important features of dynamics).

A class of systems which have a sublinear increase of the
information are the systems with an infinite invariant measure
(see Theorem \ref{teo:aic-mi}). An important subclass of these
consists of maps with an indifferent fixed point, being an
important example the map on the interval $[0,1]$ given by
\begin{equation} \label{manne001}
x\mapsto x+x^{z}\ (\mbox{mod }1)\ \mbox{ for } z\ge 2
\end{equation}
In the above family of maps the origin is a neutrally unstable
fixed point, hence an orbit that is sent near the origin can be
trapped near the origin for long times. The resulting behaviour is
an alternation of chaotic (when the orbit stays far from the
origin the map is similar to the baker's map) and regular (when
the orbit is trapped near the origin) phases. The expected
trapping times can be modulated by varying $z$. For these reasons
this map was introduced in the physical literature in \cite{mann}
as a model of intermittent turbulent behaviour in fluid dynamics,
and the particular statistical properties of the orbits of these
maps were used in different fields to model intermittent phenomena
(see for example \cite{Grigolini}).

From the information point of view these maps exhibit a behaviour
that is between the fully chaotic (positive entropy) and the
regular one (the dynamics is predictable, the information needed
to describe it increases with time $n$ at most as $\log (n)$).
This was first discovered in \cite{GW} in a piecewise linear
example (see Section \ref{sec:plm}). Their seminal, short paper
however does not present a complete mathematical proof of this
fact (a lower bound, like Theorem \ref{teo:finale2} is not
proved). Some further study was made in \cite{gal1} (where the
lower bound was proved under some assumption on $z$) and in
\cite{bo-tesi}.

This paper considers both $C^1$ and piecewise linear (PL) classes
of interval maps with an indifferent fixed point. The $C^1$ case
is studied using techniques that are different from the techniques
used in the previous literature.

The main result of this paper implies that, in mean with respect
to any absolutely continuous probability measure the information
of the "Manne\-vil\-le-Pomeau like"(see Definition \ref{def:mp})
class of maps, for which equation (\ref{manne001}) is an example,
behaves for $z>2$ like $$\E[I(x,n,Z)] \sim n^\frac{1}{z-1}$$ that
is as a power law with exponent less than 1 (see Theorem
\ref{risultato}). This shows in particular that the behaviour of
the information content for these maps depends only on the local
behaviour of the map near the neutrally unstable fixed point at
the origin. Some results about the pointwise behaviour of
$I(x,n,Z)$ are also given (see Proposition \ref{prop:pw}).

Moreover we study a class of piecewise linear maps, extending
the results of \cite{GW}, finding in particular different
behaviours for the mean of the information content (see Section
\ref{sec:plm}).

All our results are based on the definitions of Section
\ref{sec:info}, where we introduce in an informal way the
\emph{Algorithmic Information Content} of a string and the few
related facts about algorithmic information theory that we need in
the following. Then using these concepts we define the
\emph{local} and \emph{global chaos indexes}, which measure the
local and average power law behaviour of the information
$I(x,n,Z)$ (roughly speaking when $I(x,n,Z)\sim n^{\alpha }$ then
the index is $\alpha$). The invariance properties of these indexes
and the results on the Manneville-Pomeau maps imply, as a simple
corollary, that two Manneville-Pomeau maps with different
parameter $z$ cannot be absolutely continuously conjugate
(Corollary \ref{coniugio}).

\section{Information measures and chaotic \\ dynamical systems} \label{sec:info}

The method we use to study a chaotic dynamical system is based on
the idea of a measure of the information contained in the orbits
of the system.

Let $\mathcal A$ be a finite alphabet and let ${\mathcal A}^n$ be
the set of all strings of length $n$ written with letters from
$\mathcal A$. Then by ${\mathcal A}^*$ we denote the set of finite
strings of any length, that is
$${\mathcal A}^* = \bigcup_{n=0}^\infty {\mathcal A}^n$$
Given a string $s\in {\mathcal A}^*$, the intuitive idea of
information contained in $s$ is the length of the smallest binary
message from which it is possible to reconstruct $s$. Thus,
formally, the information $I$ is a function $ I:\mathcal{A}^* \to
\N.$

One of the most important measures for the information content is
the \emph{ Algorithmic Information Content (AIC)}. In order to
define it, it is necessary to define the notion of partial
recursive function. We limit ourselves to give an intuitive idea
which is very close to the formal definition. We can consider a
partial recursive function as a computer $C$ which takes a program
$p$ (namely a binary string) as an input, performs some
computations, and gives a string $s=C(p)$, written on the given
alphabet $\mathcal{A}$, as an output. The $AIC$ of a string $s$ is
defined as the shortest binary program $p$ which gives $s$ as its
output, namely
\begin{equation} \label{eq:aic}
{AIC}(s,C)=\min \{|p|:C(p)=s\}
\end{equation}
where $|\cdot|$ denotes the binary length of the program $p$. Up
to now the Algorithmic Information Content depends on $C$, but
there is a class of computing machines that allows a definition of
information content independent on the particular computer up to a
constant. We require that our computer is a universal computing
machine. Roughly speaking, a computing machine is called
\emph{universal} if it can simulate any other machine if
appropriately programmed. That is, $C$ is universal if for each
other computer $D$ there is a program $p_{CD}$ such that for each
program  $p$ it holds $D(p)=C(p_{CD}p)$. In particular the
computers we use every day are universal computing machines,
provided that we assume that they have virtually infinite memory.
For a precise definition see for example \cite{livi} or
\cite{cha}. We have the following theorem

\begin{teorema}[\cite{kolmogorov}] \label{teo:aic-kolm}
If $C$ and $C^{\prime }$ are universal computing machines then
$$\left| {AIC}(s,C)-{AIC}(s,C^{\prime })\right| \leq K\left(
C,C^{\prime }\right)$$ where $K\left( C,C^{\prime }\right)$ is a
constant which depends only on $C$ and $C^{\prime }$.
\end{teorema}

This theorem implies that the information content ${AIC}$ of $s$
with respect to $C$ depends only on $s$ up to a fixed constant,
then its asymptotic behaviour does not depend on the choice of
$C$. For this reason from now on we will write ${AIC}(s)$ instead
of ${AIC}(s,C)$.

The shortest program which gives a string as its output is a sort
of encoding of the string, and the information which is necessary
to reconstruct the string is contained in the program. From this
point of view the computer can also be seen as a decoder.
Unfortunately the coding procedure associated to the Algorithmic
Information Content cannot be performed by any algorithm. This is
a very deep statement and, in some sense, it is equivalent to the
Turing halting problem or to the G\"odel incompleteness theorem.
Then the Algorithmic Information Content is a function not
computable by any algorithm. Hence in computations one tries to
approximate from above the $AIC$ of a string by means of some
algorithm. This leads to consider the theory of \emph{compression
algorithms}, algorithms that encode an original string $s\in \A^*$
into a compressed version of it $M(s)\in \{ 0,1 \}^*$ in a
reversible way (i.e., there is an  other algorithm that is able to
recover the original string from the coded string). Using a
compression algorithm $M$ one defines the information content of
the string $s$ with respect to the encoding procedure $M$ as
$I_M(s)=|M(s)|$, where $|\cdot|$ denotes the binary length of
$M(s)$. We remark that the Algorithmic Information Content of a
string is up to a constant less than or equal to the information
content as it is computed by some compression algorithm. That is,
for each compression algorithm $M$  there is a constant $C_M$ such
that for each $s$ it holds $I_M(s)\geq AIC(s) +C_M$. This is
because each universal computing machine can be programmed also to
perform any coding-decoding technique, and the length of this
program represents $C_M$.

To study a chaotic dynamical system we have to consider the
asymptotic behaviour of its orbits, hence if we want to consider
in some sense the information contained in the orbits of the
system, we have to deal with infinite strings. Then let
$\Omega:={\mathcal A}^\N$ denote the set of infinite strings
$\omega$ with letters from ${\mathcal A}$.  To study the
information contained in a infinite string $\omega$, we study what
is the asymptotic behaviour of the information contained in the
first $n$ symbols of $\omega$ as $n$ increases. Hence we study the
asymptotic behaviour of the function
$AIC(\omega,n):=AIC(\omega^n)$ as $n$ increases, where
$\omega^n=(\omega_0,\dots,\omega_{n-1})$ is the string given by
the first $n$ symbols of $\omega$. By the methods of symbolic
dynamics we will associate an infinite string to an orbit and then
the above idea will be applied to dynamical systems.

\subsection{Application to dynamical systems} \label{sec:dynsyst}

Let $(X,\B,T)$ be a dynamical system. $X$ is assumed to be a
compact metric space, $\B$ is the Borel $\sigma$-algebra, and $T$
is a $\B$-measurable map from $X$ to itself. Let $\mu$ be a
$T$-invariant measure on $(X,\B)$. We do not suppose that
$\mu(X)<\infty$, but we always assume $\mu$ to be $\sigma$-finite
and conservative\footnote{A system is conservative if the set of
wandering points has zero measure (see \cite{aaronson}).}.

Let $Z=\{I_0,\dots,I_{N-1}\}$ be a finite measurable partition of
$X$, and let $\A=\{0,\dots,N-1\}$ be the associated finite
alphabet. Then the symbolic representation of the system
$(X,T,\mu)$ is given by the function $\varphi_Z : X\to
\Omega=\A^\N$ associating to each point a symbolic orbit with
respect to $Z$, defined by
\begin{equation} \label{eq:symb}
\varphi_Z(x)=\omega \ \Longleftrightarrow \ T^i(x) \in
I_{\omega_i} \quad \forall \ i\in \N
\end{equation}

The image $\varphi_Z(X)$ is a subset of $\Omega$ which is
invariant under the usual shift map $\tau$ on $\Omega$. The
function $\varphi_Z$ induces on $\varphi_Z(X)$ a measure
corresponding to $\mu$. Of course similar symbolic representations
can also be defined for countable partitions (with a countable
symbolic alphabet).

At this point it is possible to apply the notion of information to
the orbits of the system. We obtain quantities dependent on a
given partition $Z$ of $X$. The following
definition and the following theorem, given in \cite{brudno}, shows the relation between AIC
and entropy.

\begin{definizione} \label{def:complpts} The \emph{information content
$AIC(x,n,Z)$} of $n$ steps of the orbit of $x$ with respect to $Z$
is defined as  $AIC(x,n,Z)= AIC((\varphi_Z(x))^n)$. Analogously,
the {\it complexity} $K(x,Z)$ of a point $x\in X$ with respect to
$Z$ is given by
\begin{equation} \label{eq:complpts}
K(x,Z)=\limsup\limits_{n\to \infty} \ \frac{AIC(x,n,Z)}{n}.
\end{equation}
\end{definizione}

The complexity of a sequence is related to the Shannon entropy of
the information source that has produced the sequence. Hence in
the theory of dynamical systems it is possible to relate the
complexity to the Kolmogorov-Sinai entropy of the system. The
following theorem holds

\begin{teorema}[\cite{brudno},\cite{white}] \label{teo:ks-aic}
Let $(X,T)$ be a dynamical system and $\mu$ a $T$-invariant
probability measure on $X$. Given a finite measurable partition
$Z$ of $X$, it holds $$\int_X \ K(x,Z) \ d\mu(x) = h_\mu (T,Z)$$
where $h_\mu (T,Z)$ denotes the Kolmogorov-Sinai entropy of the
system relative to the partition $Z$. If the dynamical system
$(X,T,\mu)$ is ergodic then for $\mu$-almost all $x\in X$
\begin{equation} \label{W}
K(x,Z) = \liminf\limits_{n\to\infty} \frac{AIC(x,n,Z)}{n} = h_\mu
(T,Z)
\end{equation}
\end{teorema}

In systems with an infinite measure we have a behaviour of the
information that is similar to zero entropy systems, indeed the
following holds.

\begin{teorema}[\cite{bo-tesi}] \label{teo:aic-mi}
Let $(X,T)$ be a dynamical system and $\mu$ a $T$-invariant
infinite ergodic measure. Given a finite measurable partition $Z$
of $X$, for $\mu$-almost all $x\in X$ it holds $K(x,Z)=0$.
\end{teorema}

These theorems tell us what we can expect for the asymptotic
behaviour of the information content of a typical orbit of an
ergodic dynamical system. If the system has an invariant
probability measure $\mu$ with positive Kolmogorov-Sinai entropy
relative to a partition $Z$, then for $\mu$-almost all $x\in X$ it
holds\footnote{{\sc Notations:} Here and in the sequel, for two
sequences $a_n$ and $b_n$ we shall write $a_n \sim b_n$ if the
quotient $a_n/b_n$ tends to unity as $n\to \infty$. Moreover, the
notation $a_n\approx b_n$ means that $a_n/b_n=O(1)$ as well as
$b_n/a_n=O(1)$ for $n\to \infty$. Moreover we shall write
$a_n\preceq b_n $ if $a_n=O(b_n)$.}
\begin{equation} \label{eq:aic-asint}
AIC(x,n,Z)\sim n \ h_\mu (T,Z)
\end{equation} If instead the
ergodic measure $\mu$ has null Kolmogorov-Sinai entropy or it is
an infinite measure, then for $\mu$-almost all $x\in X$
\begin{equation} \label{eq:aic-asint2}
AIC(x,n,Z)=o(n)
\end{equation}
for any finite partition $Z$. In this second situation we
introduce the notions of \emph{local} and \emph{global chaos
indexes} to classify dynamical systems according to the asymptotic
behaviour of $AIC(x,n,Z)$.

Let $(X,T)$ be a dynamical system and $\mu$ a $T$-invariant
measure. Let $\nu$ be a probability measure on $X$ equivalent to
$\mu$ (each one is absolutely continuous with respect to the
other). We refer to $\nu$ as a "reference measure". Of course when
$\mu$ is a probability measure, we can set $\nu=\mu$. Let $Z$ be a
finite partition of $X$.

\begin{definizione} \label{def:global-ind}
The \emph{upper chaos index $\oq (T,Z,\nu)$ with respect to $Z$}
is given by
$$\oq(T,Z,\nu)= \inf \Big\{ q\in (0,1) \ /\ \limsup\limits_{n\to
\infty} \ \int_X \ \frac{AIC(x,n,Z)}{n^q} \ d\nu(x) = 0 \Big\} =$$
$$= \sup \Big\{ q\in (0,1) \ /\ \limsup\limits_{n\to \infty} \ \int_X
\ \frac{AIC(x,n,Z)}{n^q} \ d\nu(\omega) = \infty \Big\}$$ In the
same way the \emph{lower chaos index $\uq(T,Z,\nu)$ with respect
to $Z$} is defined using the inferior limit instead of the
superior limit.
\end{definizione}

We remark that in principle the chaos index may depend on the
choice of $\nu$. In principle, even if we consider equivalent
measures, the index may change. However we will see that in the
class of map we are interested to study this does not happen.
Moreover we remark that in the examples we study there is a
natural choice of the reference measure, the Lebesgue measure.

We  also define the {\it upper} and {\it lower local chaos
indexes}.

\begin{definizione} \label{def:local-ind}
The {\it upper local chaos index $\oq(T,x,Z)$} is defined as
$$\oq(T,x,Z)= \inf \Big\{ q\in (0,1) \ /\ \limsup\limits_{n\to
\infty} \ \frac{AIC(x,n,Z)}{n^q} = 0 \Big\} = $$
$$= \sup \Big\{ q\in (0,1) \ /\ \limsup\limits_{n\to \infty} \
\frac{AIC(x,n,Z)}{n^q} = \infty \Big\}$$ In the same way the {\it
lower local chaos index $\uq(T,x,Z)$} is defined using the
inferior limit instead of the superior limit.
\end{definizione}

\begin{teorema}[\cite{bo-tesi}] \label{teo:indici}
Let $\mu$ be ergodic and invariant for the dynamical system
$(X,T)$ and let $\nu$ be a reference measure. For any finite
partition $Z$ of $X$, the local indexes are a.e.-constant, that is
for $\mu$-almost all $x\in X$ it holds $\oq(T,x,Z)=\oq(T,Z)$ and
$\uq(T,x,Z)=\uq(T,Z)$. Moreover
$$\uq(T,Z) \le \uq(T,Z,\nu) \le \oq(T,Z,\nu)$$
\end{teorema}

Now we want to get rid of the dependence on the partition $Z$ and
define an indicator that is independent on the partition. We will
state a definition for maps over the interval. Further
generalizations are possible (see \cite{licatone}) but for the
sake of simplicity here we will restrict to interval maps.

To avoid pathologies coming from very complicated partitions (see
the "negative results" in \cite{blume}) we will consider a class
of admissible partitions and take the supremum over this class.

\begin{definizione}\label{XXX}
A partition $Z=\{I_0,...,I_{N-1} \}$ of $[0,1]$ is called
\emph{admissible} if it is made of a finite set of intervals, i.e.
each $I_i$ is an interval. The \emph{upper} and \emph{lower global
weak chaos indexes} of $([0,1],T,\nu)$ are defined by
$$\oq(T,\nu)= \sup\limits_{Z\in \ adm.\ part.}\oq(T,Z,\nu)$$
$$\uq(T,\nu)= \sup\limits_{Z\in \ adm.\ part.}\uq(T,Z,\nu)$$
\end{definizione}

What we obtained is a weak chaos index for maps of the interval
that is in general invariant for a bi-Lipschitz conjugacy. Let us
consider the Lebesgue measure $m$ on the unit interval.

\begin{teorema}
If $([0,1],T)$ and $([0,1],T')$ are conjugated by a bi-Lipschitz
homeomorphism $\pi$, then $\oq(T,m)=\oq(T',m),
\uq(T,m)=\uq(T',m)$.
\end{teorema}

\noindent {\bf Proof.} It is clear that a homeomorphism sends an
admissible partition $Z$ to an admissible partition $\pi(Z)$. We
have that $$\int_{[0,1]} \ \frac{AIC(\pi(x),n,\pi(Z) )}{n^q} \ dx=\int_{[0,1]} \
\frac{AIC(x,n,Z)}{n^q} \ f(x)\ dx $$ where $f(x)\in L^\infty
[0,1]$ and then we can estimate one index in function of the
other. \qed

\section{The ``Manneville-Pomeau like'' maps} \label{sec:mp}

We apply our techniques, based on the asymptotic behaviour of the
information content of symbolic orbits, to a family of interval
differentiable maps $T:[0,1]\to [0,1]$ with an indifferent fixed
point.

\begin{figure}[h]
\begin{tabular}{ll}
{{ \psfig{figure=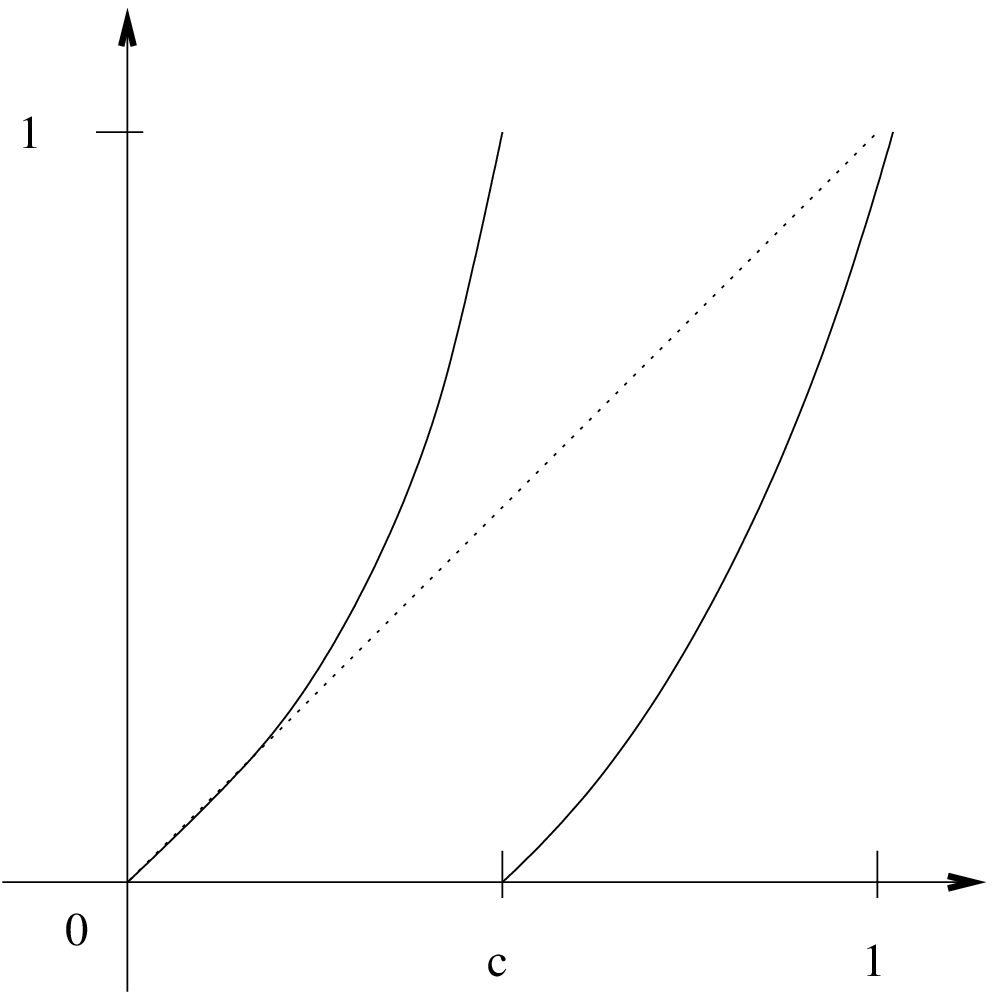,width=5cm,angle=0}}} & {
\raggedleft{\psfig{figure=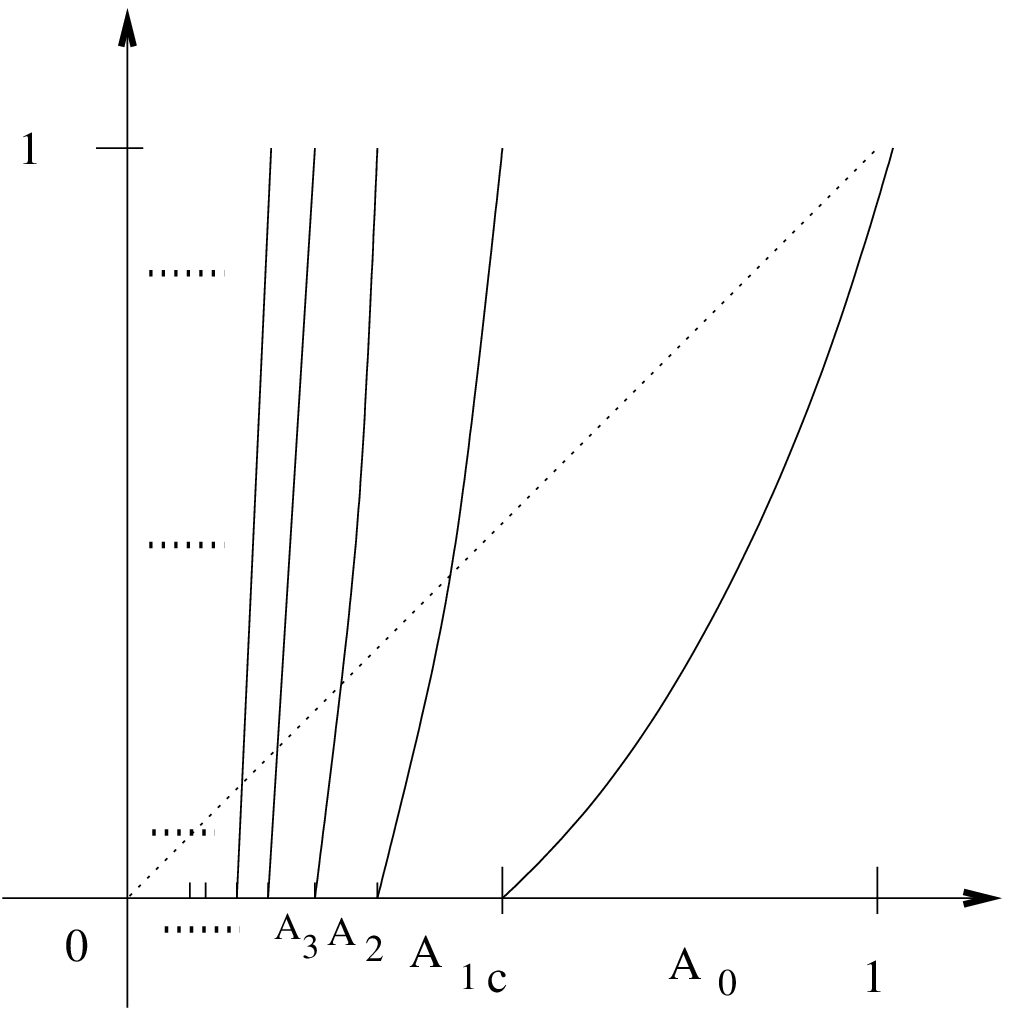,width=5cm,angle=0}}}
\end{tabular}
\caption{\it An  example of the  graph of an MP map and its
induced map.} \label{figmann}
\end{figure}

\begin{definizione} \label{def:mp}
We say that a map $T:[0,1]\to [0,1]$ is a \emph{Manneville-Pomeau
map (MP map)} with exponent $z$  if it satisfies the following
conditions:
\begin{enumerate}
    \item there is $c\in (0,1)$ such that, if $I_0=[0,c]$ and
    $I_1=(c,1]$, then $T\big|_{(0,c)}$ and $T\big|_{(c,1)}$ extend
    to $C^1$ diffeomorphisms, $T(I_0)=[0,1]$, $T(I_1)=(0,1]$ and $T(0)=0$;

    \item there is $\lambda >1$ such that $T'\ge \lambda$ on
    $I_1$, whereas $T'>1$ on $(0,c]$ and $T'(0)=1$;

    \item the map $T$ has the following behaviour when $x\to 0^+$
    $$T(x)=x+rx^z +o(x^z)$$
    for some constant $r>0$ and $z>1$ (see left part of Figure \ref{figmann}).
\end{enumerate}
\end{definizione}

This family of maps is well known and many statistical (the decay
of correlations, the central limit theorem and the phenomenon of
phase transitions) and ergodic (exactness, rational ergodicity,
mixing and the return time sequences) properties have been deeply
analyzed. Most of these studies are made for $z< 2$, on the
contrary we are mostly interested in the case $z\geq 2$.

For $z\in (1,2)$ there is a unique absolutely continuous (with
respect to Lebesgue measure) probability measure $\mu$ that is
$T$-invariant, moreover $\mu$ is a Sinai-Ruelle-Bowen measure, it
is exact and its Kolmogorov-Sinai entropy satisfies
$$h_\mu (T) = \int_0^1 \log |T'(x)| \ d\mu(x) > 0$$
From the statistical point of view, it is known that some changes
happen when $z$ crosses the value $z=\frac{3}{2}$. However from
our point of view these changes are not relevant.

Applying Theorem \ref{teo:ks-aic} and the consequent equation
(\ref{eq:aic-asint}) we obtain that if $T$ is a Manneville-Pomeau
map with $z\in (1,2)$, and $Z$ is the generating partition $Z=\{
I_0,I_1 \}$, then for the absolutely continuous $T$-invariant
probability measure $\mu$ it holds
$$AIC(x,n,Z)\sim n\ h_\mu (T)$$
for $\mu$-almost all $x\in [0,1]$.

Much more delicate is to study the case $z\ge 2$. Indeed for these
values of the parameter the only absolutely continuous
$T$-invariant measure $\mu$ is infinite. Hence in this case we
obtain from Theorem \ref{teo:aic-mi} and equation
(\ref{eq:aic-asint2}) that $AIC(x,n,Z)=o(n)$ for $\mu$-almost all
$x\in [0,1]$, for any finite partition $Z$. To classify these maps
we study the behaviour of the chaos indexes introduced in Section
\ref{sec:info}.

The ergodic properties of the infinite measure $\mu$ for MP maps
with $z\ge 2$ have been studied in \cite{thaler2} and
\cite{thaler}, applying the theory of infinite ergodic measures
(see \cite{aaronson}). In particular the measure $\mu$ is shown to
be exact and rationally ergodic, with estimates for the return
time sequences. Putting together the results of \cite{thaler} and
the Darling-Kac Theorem (\cite{aaronson}) we obtain

\begin{teorema}[Aaronson-Darling-Kac-Thaler] \label{teo:adkt}
Let $T$ be a Manneville-Pomeau map with $z\ge 2$ and $\mu$ the
infinite absolutely continuous $T$-invariant measure. Then for all
Borel measurable $B\subset [0,1]$ with $\mu(B)<\infty$ it holds
$$\frac{1}{a_n} \ \sum_{j=0}^{n-1} \ \chi_B \circ T^j \
\stackrel{\mathcal{L}}{\Longrightarrow} \ \mu(B)\ Y_\alpha$$ where
the convergence is in distribution with respect to all absolutely
continuous Borel probability measures on $[0,1]$, $Y_\alpha$ is a
positive random variable distributed according to the normalized
Mittag-Leffler law of order $\alpha$ and $\alpha = \frac{1}{z-1}$.
Moreover it holds
\begin{itemize}
    \item $a_n \sim \frac{n}{\log n}$ if $z=2$;
    \item $a_n \sim n^{\frac{1}{z-1}}$ if $z>2$.
\end{itemize}
\end{teorema}

The statistical distribution of the frequencies of visits to
subsets of $[0,1]$ is the fundamental tool to obtain the behaviour
of the $AIC$ of orbits of the system.

Our plan is the following:
\begin{itemize}
    \item[(i)] to consider a symbolic representation of the system
induced by the choice of a partition;
    \item[(ii)] to estimate the information content of the orbits using a
particular encoding as a compression algorithm;
     \item[(iii)] to show that this information content has the same
 average asymptotic behaviour as the Algorithmic Information Content.
\end{itemize}

To obtain a symbolic representation of the system we will consider
a finite admissible partition (see Definition \ref{XXX}). First we
will consider the case of a partition made of two intervals
$I_0=[0,c]$ and $I_1=(c,1]$, then we will show that the general
case is similar. In this first particular case, the alphabet
associated to the partition $Z$ is ${\cal A}=\{0,1\}$.

Let us consider a MP map. As said before the presence of the
indifferent fixed point at the origin implies that a typical orbit
will spend much time near the origin, since it moves away from it
very slowly, and will spend the rest of time around in the
interval. Then at some time it will come again close to the origin
and again stay near the origin a lot of time. This repeats over
and over again. This fact implies that a typical symbolic orbit
will have a lot of symbols equal to "0" and some equal to "1".
Moreover being the derivative of the map bounded from 1 in $I_1$,
the dynamics in this interval is "fully" chaotic, and this implies
a sort of renewal process when the orbit of the point is in $I_1$
(see Section \ref{sec:plm}). The results we will see are in some
sense reminiscent of the theory of renewal processes, but this
theory fails for MP maps, since they are not isomorphic to a
Markov chain.

We now introduce the encoding we use to estimate the information.
A typical sequence $\omega =\varphi_Z(x) \in \{0,1\}^{\N}$ will
look like\footnote{This symbolic representation is well defined
for $m$-a.e. point $x\in [0,1]$ (analogously to the case of the
dyadic numbers for the Bernoulli shift).}
$$\omega=(000100000000110010000000001101\dots)$$
so it is possible to compress its first 30 symbols in the
following string
\begin{equation} \label{eq:compress}
s(\omega^{30})=(3802901)
\end{equation}
where we have just written how many "0"s there are between two
consecutive "1"s. Since the number of consecutive zeros can be as
high as we want, the string $s$ is a string with digits coming
from an infinite alphabet (each digit is a natural number). Since
we want to deal with strings coming from a finite alphabet we
codify $s$ into a binary string $s'$, simply using the standard
binary representation of natural numbers and writing it with the
usual prefix-code (see below for an example). In this way a number $n\in \N$ is encoded by a
binary string $(a_0,a_1,\dots,a_m)$ such that
$$n=2^{m+1}-1+\sum_{j=0}^m a_j 2^j$$
where $m+1 = \lfloor \log_2 (n+1) \rfloor$ (denoting $\lfloor
\cdot \rfloor$ the inferior integral part of a number), and this
binary string is written using the prefix-code given by
$$\bar n = (1^{m+1}, 0, a_0, \dots, a_m)$$
where $1^{m+1}$ denotes the symbol "1" repeated $m+1$ times. In
this way, leaving unchanged the symbols "0", we obtain for our
example
$$s'(\omega^{30})=(11000\ 1110001\ 0\ 101\ 1110010\ 100)$$
Hence $s'(\omega^n)$ is an encoding of the string $\omega^n$ (see
Section \ref{sec:info}).

The information function $I$ associated to this compression is
given by
\begin{equation} \label{eq:inf-func}
I_{s'}(\omega^n) =|s'(\omega^n)|= \sum_{j=0}^{|s(\omega^n)|-1} \
2\ \lfloor\ \log_2 (1+s(\omega^n)_j)\ \rfloor +1
\end{equation}
where $|\cdot|$ denotes the length of the string.

Let $N_n:\Omega \to \N$ be the sequence of functions defined by
\begin{equation} \label{eq:passaggi}
N_n(\omega)=\# \{ i \ | \ \omega_i\neq 0, \quad i=0,\dots,n-1 \}
\end{equation}
that is the number of passages of the orbit outside the interval
$I_0$ in the first $n$ steps. Then it is easy to realize that in
the case of a partition made by two intervals
$N_n(\omega)=|s(\omega^n)|$. In this case we now prove a stronger
relation between $I_{s'}(\omega^n)$ and $N_n(\omega)$.

\begin{lemma} \label{teo:stima-inf}
For any sequence $\omega \in \{0,1\}^{\N}$ it holds
\begin{equation*}
N_n(\omega) + 2\ \log_2 (n-N_n(\omega) +1) \leq I_{s'}(\omega^n)
\leq N_n(\omega) + 2\ N_n(\omega) \log_2 \left(
\frac{n}{N_n(\omega)} \right) \label{eq:stimainf}
\end{equation*}
up to an additive constant given by the possibilities $\omega_0=1$
and $\omega_{n-1}=0$, and by the presence of the inferior integral
part in the definition of $I_{s'}(\omega^n)$.
\end{lemma}

\noindent {\bf Proof.} It is enough to prove the lemma for
$\omega_0=0$ and $\omega_{n-1}=1$. Otherwise simply add a
constant.

In general $N_n = n-h$, for some $h < n$. The compression of such
strings is then $N_n$-symbols long. Moreover the compression is
such that $\sum_{j=0}^{N_n -1} s(\omega^n)_j = h$. We now want to
find the maximum and the minimum of the function
\[
\sum_{j=0}^{N_n -1} \ 2\ \log_2 (1+s(\omega^n)_j)
\]
with the condition $\sum_{j=0}^{N_n -1} s(\omega^n)_j = h$. The
maximum is attained for equal $s(\omega^n)_j \not= 0$, and the
minimum for all the $s(\omega^n)_j =0$ but one which is equal to
$h$. Then the maximum is given by $s(\omega^n)_j = \frac{h}{n-h}$
for all $j$, and the information content is given by
$\sum_{j=0}^{n-h-1} \ \left[ 1+2\ \log_2 \left( 1+ \frac{h}{n-h}
\right) \right] = N_n + 2\ N_n \log_2 \left( \frac{n}{N_n}
\right)$, and the minimum is given by $(n-h-1) + 1+2\ \log_2 (h+1)
= N_n +2\ \log_2 (n-N_n +1)$. Hence the lemma is proved. \qed

\vskip 0.5cm Let us see what can be done when we have some general
admissible partition made of a finite number of intervals. Let
$Z=\{ I_0,\dots,I_{N-1} \}$ be such a partition and let ${\cal A}$
be the associated finite alphabet. In this case we assume $\A$ to
be made of the symbol "0" for the interval $I_0$, and of letters
(or any other kind of symbols different from natural numbers) for
the other intervals.

We slightly modify the previous coding procedure as follows. We
have strings $\omega ^n\in {\cal A}^n$, and define
$s(\omega^n)=(x_1,x_2,\dots,x_m)$, where $x_i\in \N \cup {\cal
A}$. In the sense that we codify as numbers the occurrences of the
symbol "0" as before, in a way that if at place $i$ there are
$n_i$ consecutive "0"s we write the number $n_i$. The other
symbols are left unchanged. For example, if ${\cal A}=\{0,A,B\}$
and $\omega^n=(00ABB000A)$ then $s(\omega^n )=(2ABB3A)$. Then we
define $s'$ as before, by the standard binary encoding of the
natural numbers, obtaining a string $s'$ written in the alphabet
$\{ 0,1,A,B \}$. In the previous example, we obtain
$s'(\omega^n)=(101\ ABB\ 11000\ A)$.

Then we can easily estimate the information function $I(\omega^n)$
in this case, by noting that the only difference with the previous
case of a partition with only two intervals is that now when the
symbol in $\omega^n$ is different from "0" we have to explicitly
specify it, so that in $s'(\omega^n)$ the symbols "0" are replaced
by the explicit strings, that is "$ABB$" and "$A$" in our example.
Since there are $N_n(\omega)$ such symbols different from "0" or
"1" in $s'(\omega^n)$, we need at most $N_n(\omega) \log_2 (\# \A
-1)$ bits more than in the previous case with a partition with
only two intervals. Note that in the previous case $\log_2 (\# \A
-1) =0$, hence this new term vanishes. We can now state the result

\begin{lemma}\label{3.5}
Let $\omega^n\in {\cal A}^n$ and $s'$ be the encoding as above,
then it holds
$$I_{s'}(\omega ^n)\leq\ N_n(\omega) + 2\ N_n(\omega) \log_2 \left(
\frac{n}{N_n(\omega)} \right) + N_n(\omega) \log_2 (\# \A -1)$$
where $\#{\cal A}$ is the cardinality of ${\cal A}$.
\end{lemma}

This lemma implies that the behaviour of the information content
is given by the asymptotic behaviour of the functions $N_n$. To
estimate the behaviour of $N_n$ first of all notice that if $x$ is
a point in $[0,1]$ such that $\varphi_Z(x)=\omega$, then
$$N_n(\omega)= \sum_{j=0}^{n-1} \ \chi_{_{[0,1]\setminus I_0}} (T^j(x))$$
Then we can apply Theorem \ref{teo:adkt} to the sequence $N_n$
with respect to any measure $\nu$ on $\Omega$ induced by an
absolutely continuous probability measure on $[0,1]$. This gives
the estimates
\begin{equation} \label{eq:asint-n}
\E_\nu [N_n(\omega)] \sim \left\{
\begin{array}{ll}
    \frac{n}{\log n}, & \hbox{z=2;} \\[0.2cm]
    n^{\frac{1}{z-1}}, & \hbox{z $>$ 2.} \\
\end{array}
\right.
\end{equation}
obtained by the asymptotic behaviour of the sequence $a_n$ in
Theorem \ref{teo:adkt}, where $\E_\nu [\cdot]$ denotes the mean
with respect to the measure $\nu$. By the above estimates we have

\begin{prop} \label{teo:finale1}
Let $T$ be a Manneville-Pomeau map with parameter $z\ge 2$. Let
$Z$ be an admissible finite partition. Then for any probability
measure $\nu$ on $\Omega$ induced by an absolutely continuous
probability measure on $[0,1]$ through the symbolic representation
$\varphi_Z$, it holds
$$\begin{array}{ll}
 \E_\nu [AIC(\omega^n)]\preceq   \E_\nu [I(\omega^n)] \preceq \ n, & \hbox{z=2} \\[0.3cm]
 \E_\nu [AIC(\omega^n)]\preceq    \E_\nu [I(\omega^n)] \preceq \ n^\frac{1}{z-1} \log (n), & \hbox{z $>$ 2}. \\
\end{array}
$$
\end{prop}

Proposition \ref{teo:finale1} gives an estimate from above for the
asymptotic behaviour of $\E_\nu [AIC(\omega^n)]$. Now we want to
calculate a lower estimate. For this we consider a map that is
derived by the MP map (it is an induced map) and which gives a
symbolic dynamics that is similar to the compressed string
$s(\omega^n)$. Proving that such a map has positive entropy we
prove (using Theorem \ref{teo:ks-aic}) that no further drastic
compression is possible.

Let us consider again the partition $Z=\{[0,c],(c,1]\}$, we will
give a lower estimate to $\underline q(T,Z,\nu)$. The induced
version of a MP map is obtained studying the passages through
$I_1$. For any $x\in [0,1]$ let $\tau(x)$ be the \emph{time of the
first passage through $I_1$}, that is
\begin{equation} \label{eq:first-pass}
\tau(x) = 1+ \min \ \left\{ n\ge 0 \ |\ T^n(x)\in I_1 \right\}
\end{equation}
The level sets of the function $\tau(x)$ are a partition of
$(0,1]$ into intervals $(A_n)_{n\ge 1}$ defined as
\begin{equation} \label{eq:level}
A_n = \left\{ x\in [0,1] \ |\ \tau(x)=n \right\}
\end{equation}
Using $\tau(x)$ we define the \emph{induced map} $G:[0,1] \to
[0,1]$ by
\begin{equation} \label{eq:ind-map}
G(x)=T^{\tau(x)}(x)
\end{equation}
hence $G|_{A_n} (x) = T^n(x)$ and $G(A_n)=(0,1]$, finally we
define $G(0)=0$ (see right part of Figure \ref{figmann}).

By the definition of the induced map, we have the following
proposition.

\begin{prop} \label{prop:1}
Let $T$ be a Manneville-Pomeau map and $G$ its induced version on
$I_1$. Let $Z$ be the finite partition $Z=\{ I_0,I_1 \}=\{
[0,c],(c,1] \}$ and $P$ be the countable partition of $[0,1]$ into
the intervals $(A_n)_{n\ge 0}$, with $A_0=\{[c,1] \}$ and $A_n$
for $n\ge 1$ defined as in equation (\ref{eq:level}). Then for a
point $x$ it holds $\varphi_Z(x)=\omega$ if and only if
$\varphi_P(x)=s(\omega)$\footnote{By the definition of a
compression algorithm, $s(\omega)$ makes sense only if $\omega$ is
a finite string, hence this equation has to be intended valid for
any finite substring $\omega^n$ of the sequence $\omega$.}, where
$\varphi_P$ denotes the symbolic representation of the induced map
$G$ for the partition $P$.
\end{prop}

The entropy of the induced map is positive. This is obtained
putting together the results of Section 10 in \cite{isola} and
Theorem 22 in \cite{walters}.

\begin{teorema}[Isola-Walters] \label{teo:iw}
Let $T$ be a Manneville-Pomeau map of parameter $z$ and $G$ be its
induced version. Then there exists an absolutely continuous
$G$-invariant probability measure $\rho$ on $[0,1]$ such that: (i)
$\rho$ is exact for $G$; (ii) $\rho (A_n) \sim
n^{-1-\frac{1}{z-1}}$, hence the entropy of the partition $P$
respect to $\rho$ is finite; (iii) $\infty > h_\rho (G) = \int_0^1
\ \log |G'| \ d\rho(x) \ > 0$.
\end{teorema}

\begin{prop} \label{teo:finale2}
Let $T$ be a Manneville-Pomeau map with exponent $z\ge 2$ and let
$Z$ be the finite partition $\{ [0,c], (c,1] \}$. Then using as
reference measure any absolutely continuous probability measure
$\nu$ on $[0,1]$, it holds for the lower global weak chaos index
$$\uq(T,\nu)\geq \uq(T,Z,\nu) \geq \left\{
\begin{array}{ll}
    1, & \hbox{z=2} \\[0.2cm]
    \frac{1}{z-1}, & \hbox{z$>$2} \\
\end{array}
\right .$$
\end{prop}

\noindent {\bf Proof.} The inequality $\uq(T,\nu)\geq
\uq(T,Z,\nu)$ follows trivially by definition. By Proposition
\ref{prop:1} we have that if $\omega ^n$ is a symbolic orbit of
$x$ with respect to $Z$ then a symbolic orbit of the induced map
$G$  with respect to $P$ is $s(\omega ^n)$ and $|s(\omega
^n)|=N_n(\omega)$. Given a string $s=(s_1s_2...s_m)$ in $\N^*$
(like $s(\omega ^n)$) let us consider the string
$trunc_k(s)=(min(k,s_1)min(k,s_2)...min(k,s_n))$ where each digit
$s_i$ in replaced by $min(k,s_i)$. Up to a constant, for each $k$
$AIC(\omega ^n)\geq AIC(trunc_k(s(\omega ^n))$ because the string
$trunc_k(s(\omega ^n))$ can be obtained easily from $\omega ^n$ by
an algorithm.

We know that the entropy of the induced map is finite and
positive, then, since the infinite partition $P=\{A_0,A_1...\}$ is
generating, its entropy $h_\rho(G,P)$ is finite and positive
(Theorem \ref{teo:iw} (iii)), then there is a partition $P'_k$ of
the form $P'_k=\{A_0,...,A_k,\cup_{i>k} A_i\}$ such that
$h_\rho(G,P)\ge h_\rho(G,P'_k)>0$. This is because for each
absolutely continuous measure the sequence of partitions $P'_k$
converges to $P$ in the Rokhlin metric (and the entropy is
continuous with respect to change of partition in this metric).

We remark that $trunc_k(s(\omega ^n))$ is the symbolic orbit of
$x$ for the induced map, with respect to the partition $P'_k$ and
the length of this string is $N_n(\omega)$. Then by equation
(\ref{W}) of Theorem \ref{teo:ks-aic} and by Theorem
\ref{teo:adkt}, we have that given any sequence $\{ b_n \}$ of
positive numbers such that $b_n = o(a_n)$, where $a_n=\left\{
\begin{array}{ll}
    \frac {n}{log(n)}, & \hbox{z=2} \\[0.2cm]
   n^{\frac{1}{z-1}}  , & \hbox{z$>$2} \\
\end{array}
\right.$, and given any strictly positive constants $h,C,\epsilon$
such that $h<h_\rho (G,P'_k)$, there exists an integer $\bar n$
such that for all $n\ge \bar n$ it holds $\rho(S(n))>1-\epsilon$,
where
$$S(n) = \left\{ \omega \ / \ AIC(\omega^n)\ge
AIC(trunc_k(s(\omega^n)))\ge h N_n(\omega)\ ;\ N_n(\omega)\ge C
b_n \right\}$$ Hence, being $\rho$ absolutely continuous, for any
absolutely continuous probability measure $\nu$ and for all
$\delta >0$ there exists a $\bar n$ (different from above) such
that for each $n>\overline n$
$$\int_\Omega \ \frac{AIC(\omega^n)}{b_n} \ d\nu \ \ge \ hC \
\nu(S(n))> hC(1-\delta)$$ where $h,C$ are as above. Hence for all
sequences $b_n$ such that $b_n =o(a_n)$ it holds
$$\liminf\limits_{n\to \infty} \int_\Omega \ \frac{AIC(\omega^n)}{b_n} \ d\nu
=\infty$$ hence the thesis follows from the definition of $a_n$.
\qed

\vskip 0.5cm We are finally ready to prove the main result of the
paper.

\begin{teorema} \label{risultato}
If $\nu$ is an absolutely continuous probability measure, the
global weak chaos indexes of a MP map with exponent $z\ge 2$ are
given by
$$\uq(T,\nu)=\oq(T,\nu) =\left\{
\begin{array}{ll}
    1, & \hbox{z=2} \\[0.2cm]
    \frac{1}{z-1}, & \hbox{z$>$2} \\
\end{array}
\right.$$
\end{teorema}

\noindent {\bf Proof.} From Proposition \ref{teo:finale1} we have
the upper estimate
$$\oq(T,Z,\nu)\le \left\{
\begin{array}{ll}
    1, & \hbox{z=2} \\[0.2cm]
    \frac{1}{z-1}, & \hbox{z$>$2} \\
\end{array}
\right.$$ for each admissible partition. Moreover from Proposition
\ref{teo:finale2} we have the lower estimate for $\uq(T,\nu)$.\qed

\vskip 0.5cm It is well known that any two MP maps are
topologically conjugated, the following shows that this conjugacy
cannot be absolutely continuous if the exponents are different.

\begin{cor} \label{coniugio}
If $T_z$ and $T_{z'}$ are MP maps with different exponents $z\neq
z'$ then there are not absolutely continuous conjugacies between
$T_z$ and $T_{z'}$.
\end{cor}

\noindent {\bf Proof.} As said before an homeomorphism sends an
admissible partition $Z$ to an admissible partition $\pi(Z)$.
Considering the Lebesgue measure $dx$, we have that $$\int_{[0,1]}
\ \frac{AIC(\pi(x),n,\pi(Z) )}{n^q} \ dx=\int_{[0,1]} \
\frac{AIC(x,n,Z)}{n^q} \ \pi(dx)$$ \noindent where $ \pi(dx)$ is absolutely continuous. By Theorem \ref{risultato}, the
global chaos index of the MP maps is the same for each absolutely
continuous reference measure. By this and the above equation the
chaos index of the systems $([0,1],T_z,dx)$, $([0,1],T_z,\pi(dx))$
and $([0,1],T_{z'},dx)$ should be the same. Since the chaos index
depends only on $z$, and $z\neq z'$ this leads to a contradiction.
\qed

\vskip 0.5cm We now conclude our exposition of results about the
asymptotic behaviour of the Algorithmic Information Content of the
MP maps with a result about the local behaviour. From Theorem
\ref{teo:adkt} we can deduce that for almost each symbolic orbit
$\omega$ of a MP map with respect to the partition $Z=\{
[0,c],(c,1]\}$ the following relations hold
\begin{equation} \label{liminf1}
\limsup_{n\rightarrow \infty} \ \frac
{N_n(\omega)}{n^{\frac{1}{z-1}}} >0 \ \mbox{ for }\ z>2,
\end{equation}

\begin{equation}\label{liminf2}
\limsup_{n\rightarrow \infty}\ \frac {N_n(\omega)}{n \log(n)}>0 \
\mbox{ for }\ z=2.
\end{equation}

Hence applying the same techniques of the proof of Proposition
\ref{teo:finale2} to the upper local chaos index, and recalling
Theorem \ref{teo:indici} we have

\begin{prop} \label{prop:pw} If $Z=\{[0,c],(c,1]\}$ is the generating
partition with two intervals of a MP map, then for Lebesgue almost
each $x\in [0,1]$ it holds
$$\oq(T,x,Z)=\oq(T,Z)\geq \left\{
\begin{array}{ll}
    1, & \hbox{z=2} \\[0.2cm]
    \frac{1}{z-1}, & \hbox{z$>$2} \\
\end{array}
\right.$$
\end{prop}

This result has been applied in \cite{bgi} to obtain quantitative
recurrence results  for the maps of the Manneville-Pomeau family.

\section{Piecewise linear maps} \label{sec:plm}

We now study the behaviour of the Algorithmic Information Content
for sequences obtained as symbolic representation of orbits of a
class of piecewise linear maps whose properties are similar to the
MP maps. In this case there is an isomorphism with a Markov chain
with infinite states. We study the piecewise linear maps with the
following properties.

\begin{definizione} \label{def:plm}
Let $\{ \epsilon_k \}_{k\in \N}$  be a sequence of positive real
numbers, strictly monotonically decreasing and converging towards
zero, with the property that
\begin{equation} \label{eq:propdieps}
\frac{\epsilon_{k-1}-\epsilon_k}{\epsilon_{k-2}-\epsilon_{k-1}} <
1 \hskip 0.5cm \forall \ k \in \N
\end{equation}
We consider \emph{piecewise linear (PL) maps} $L:[0,1]\to [0,1]$
defined by
\begin{equation} \label{eq:plm}
L(x)=\left\{
\begin{array}{ll}
\frac{\epsilon_{k-2}-\epsilon_{k-1}}{\epsilon_{k-1}-\epsilon_k} (x
- \epsilon_k) + \epsilon_{k-1} & \quad \epsilon_k < x \leq
\epsilon_{k-1}, \quad k \geq 1 \\[0.3cm]
\frac{x-\epsilon_0}{1-\epsilon_0} & \quad \epsilon_0 < x \leq 1 \\[0.3cm]
0 & \quad x=0
\end{array}
\right.
\end{equation}
where $\epsilon_{-1} = 1$. Clearly the properties of the maps
depend on the sequence $\{ \epsilon_k \}$, hence varying the
sequence we obtain a class of PL maps.
\end{definizione}

Using the same approach of Section \ref{sec:mp} we consider the
partition $Z=\{I_0,I_1 \}$, where $I_0=[0,\epsilon_0]$ and
$I_1=(\epsilon_0,1]$. Hence again we have the symbolic
representation $\varphi_Z:[0,1]\to \Omega=\{0,1\}^*$. The general
case of an admissible partition can be treated similarly as in
Section \ref{sec:mp}. We also use the same information function
used for the Manneville-Pomeau maps, hence the same kind of
compression (see equation (\ref{eq:compress}) and the definition
of $s'(\omega^n)$), and study the information
$I_{s'}(\omega^n)=|s'(\omega^n)|$ (see equation
(\ref{eq:inf-func})).

Again Lemma \ref{teo:stima-inf} will apply, giving
\begin{equation} \label{teo:plm1}
\E_m[AIC(I(\omega^n))] \preceq \ \E_m[I(\omega^n)] \preceq \
\E_m[N_n] \log n
\end{equation}
and then to have an upper bound to the information it is
sufficient to study the behaviour of the random variables
$N_n:\Omega \to \N$. We will see that the behaviour of $N_n$
depends from the behaviour of $(\epsilon _k)$.

To obtain this we use the theory of infinite Markov chains.
Indeed, due to their piecewise linearity, these maps are
isomorphic to a Markov chain on an infinite alphabet, the natural
numbers $\N$. The Markov chain is defined on the probability space
$[0,1]$ with the Lebesgue measure, and the chain is in the state
$i$ at time $n$ if and only if $L^n(x)\in
(\epsilon_{i-1},\epsilon_{i-2}]$.

The transition matrix will look as follows:
\begin{equation} \label{eq:stoch-mat}
\left(
\begin{array}{ccccccc}
(\epsilon_{-1}-\epsilon_0) &(\epsilon_{0}-\epsilon_1)  & (\epsilon_{1}-\epsilon_2) & \cdots &  (\epsilon_{n-2}-\epsilon_{n-1})    & \cdots   & \cdots \\
1 & 0 & 0 & \cdots & 0 & 0   & \cdots \\
0 & 1 & 0 & \cdots & 0 & 0   & \cdots \\
0 & 0 & 1 & \cdots & 0 & 0   & \cdots \\
\cdots & \cdots & \cdots &  \cdots & \cdots & \cdots & \cdots
\end{array}
\right)
\end{equation}

This construction is well-known and we refer to \cite{GW} for its
definition. To study the behaviour of this Markov chain the
classical theory of Markov chains can be applied giving the
results that follows.

Let $t_0$ be the "mean recurrence time of the passage through the
interval $I_1$", in our case
\begin{equation}
t_0 = \sum_{k=1}^{+\infty} \ k \ (\epsilon_{k-1}-\epsilon_{k-2})
\label{eq:mrt}
\end{equation}
The first result is the existence of an invariant measure for the
Markov chain associated to our dynamical system.

\begin{teorema}[\cite{chung}] \label{teo:misura-stoc} There is a measure
$\bar p$ invariant for the Markov chain. The probability of the
event $k$ is defined by $ \bar p (k) = \sum_{n=0}^{+\infty} \
(\epsilon_{n+k-1}-\epsilon_{n+k-2})$. This measure is finite if
and only if the mean recurrence time $t_0$ is finite.
\end{teorema}

The following result is obtained using the theory of recurrent
events (\cite{feller}) and of power series (\cite{titchmarsh}).

\begin{teorema} \label{teo:inf}
Let $t_0$ be as above. If $t_0<\infty$ then $\E[N_n]
\sim\frac{n}{t_0} $, if instead $t_0=\infty$, then $\E[N_n]$ is an
infinite of order less than $n$. Moreover, let  $F(x)=
\sum_{r=1}^{[x]}(\epsilon_{k-1}-\epsilon_{k-2})$. If \begin{equation}\label{leggepot}F(x) \sim 1-
A x^{- \alpha}\end{equation} as $x\to \infty$, where $A$ is a constant and
$\alpha>1$ then $t_0<\infty$.

\noindent If $F(x)$ is as above and $0<\alpha
<1$, then $t_0=\infty $ and
\begin{equation}\label{teo:feller}
\E[N_n] \sim \ \frac{\sin \alpha \pi}{A \alpha \pi}\ n^{\alpha}
\end{equation}
moreover if $b_n=o(\E[N_n])$ the set $R(n)=\{ \omega\ / \
N_n(\omega)\leq b_n\}$ is such that
\begin{equation} \label{feller3}
m(R(n))\rightarrow 0.
\end{equation}

\noindent If $F(x)\sim 1-\frac{1}{\log x}$ then $t_0=\infty $ and
\begin{equation}\label{mild}
\E[N_n] \sim \log n.
\end{equation}
\end{teorema}

\noindent {\bf Proof.} The proof of the first statements is based
on a characterization of the mean $\E[N_n]$. In \cite{feller}, it
is shown that $\E[N_n] = U_n -1$, where $U_n = \sum_{i=0}^n \
u_i$, being $u_k$ the probability of being in the state "1" at
time $k$. Then, using Theorem 1 in \cite{feller}, it follows that
\[
\lim_{n\to +\infty} \frac{\E[N_n]}{n} = \lim_{n\to +\infty}
\frac{U_n}{n} = \lim_{n\to +\infty} u_n = \ \frac{1}{t_0}
\]
Then, if $t_0<+\infty$, that is "1" is ergodic, then $\E[N_n]$ is
linear on $n$. Whereas if $t_0 = +\infty$, that is "1" is a null
state, then $\E[N_n]=o(n)$.

Equations (\ref{teo:feller}) and (\ref{feller3}) follow from
\cite{feller}, Theorems 10 and 7.

Equation (\ref{mild}) is obtained using the results in
\cite{titchmarsh}, p.242, applications of Tauberian theorems, and
repeating the argument of Theorem 10 in \cite{feller}. \qed

\vskip 0.5cm Theorem \ref{teo:inf} and equation (\ref{teo:plm1})
give, as in the previous section, an upper estimate of the
Algorithmic Information Content of the symbolic orbits of our PL
maps.

The final step is to have a lower estimate. This will be done as
before, by considering an induced map $G$ of equation
(\ref{eq:ind-map}). In this case the level sets of the first
passage time are exactly the intervals $A_n=(\epsilon_{n-1},
\epsilon_{n-2}]$. In this case the $G$-invariant probability
measure on $[0,1]$ is the Lebesgue measure $m$, hence
$m(A_n)=(\epsilon_{n-1}-\epsilon_{n-2})$ for all $n\ge 1$. Since
the induced map $G$ is now isomorphic to a stochastic process of
independent and identically distributed random variables with
values on $S=\{ 1,2,3,\dots \}$, the $AIC$ of sequences in
$\Omega$ is equivalent to the information function $I$ if and only
if the entropy of $\{ A_n \}_{n\ge 1}$ with respect to $m$ is
finite, that is if and only if
\begin{equation} \label{eq:condiz}
H:=-\sum_{n\ge 1} \ m(A_n) \ \log(m(A_n)) \ < \infty
\end{equation}

\begin{teorema} \label{teo:plm2}
Let $L$ be a PL map on $[0,1]$ and $Z$ be the finite partition $\{
I_0, I_1 \}$. Then if $m$ denotes the Lebesgue measure and the
sequence $(\epsilon_{k})$ satisfies equation (\ref{leggepot}) it
holds
$$\E_m[N_n] \preceq \E_m[AIC(x,n,Z)] \preceq \E_m[N_n] \log n$$
\end{teorema}

\noindent {\bf Proof.} The first inequality follows directly from
equation (\ref{teo:plm1}). If $(\epsilon_{k})$ satisfies equation
(\ref{leggepot}) then equation (\ref{eq:condiz}) is satisfied. The
proof of the second inequality now is very similar to the proof of
Theorem \ref{teo:finale2}, we only have to use equation
(\ref{feller3}) instead of Theorem \ref{teo:adkt}. \qed

\vskip 0.5cm We now apply the results of this section to some
particular PL maps, showing the behaviour of the global chaos
indexes and the other features in term of the behaviour of the
sequence $\epsilon_k$. In the first column we specify the
behaviour of $\epsilon_k$, the related value of $t_0$ is in the
second column, then we have respectively the expected value of
$N_n$ (by \ref{teo:inf}), the entropy $H$ of the induced map, the
expected behaviour of the AIC of the orbits (by \ref{teo:plm2}),
the common value of the upper and lower global chaos indexes. The
expected values and the chaos indexes are with respect to the
Lebesgue measure. In the last example equation (\ref{eq:condiz})
is not satisfied, nevertheless the result for the Algorithmic
Information Content follows already from the upper estimate of
equation (\ref{teo:plm1}).

\vskip 0.5cm
\begin{tabular}{|l|l|l|l|l|l|}
$\epsilon _{k}$ behaviour & $t_{0}$ & $\E[N_{n}]$ & $H$ &
$\E[AIC(\omega ^{n})]$ & $q$ \\
\hline
$\epsilon _{k}\sim \frac{1}{a^{k}} \ \ a>1$ & $<\infty $ & $\sim \frac{n}{t_{0}%
}$ &$<\infty$ &$\approx n$ & $1$ \\[0.2cm]
$\epsilon _{k}\sim \frac{1}{k^{\alpha }} \ \ \alpha >1$ & $<\infty
$ & $\sim
\frac{n}{t_{0}}$ &  $<\infty$ & $\approx n$ & $1$ \\[0.2cm]
$\epsilon _{k}\sim \frac{1}{k^{\alpha }} \ \ \alpha <1$ & $\infty
$ & $\approx n^{\alpha }$ &$<\infty$ & $n^{\alpha }\preceq
\dots \preceq n^{\alpha }\log (n)$ & $\alpha $ \\[0.2cm]
$\epsilon _{k}\sim \frac{1}{\log (k)}$ & $\infty $ & $\sim \log n$
&$\infty$ & $\preceq (\log n)^2$ & $0$
\end{tabular}

\section{Conclusions} \label{sec:concl}

In this paper we considered two classes of weakly chaotic maps of
the interval with a neutrally unstable fixed point. We calculated
the information with respect to a partition and showed that this
gives an invariant to characterize different weakly chaotic
dynamics.

This kind of behaviour for dynamical systems has been largely
studied in the last years, from many different points of view. The
importance of our approach lies in the fact that whereas the
results given here are theoretical, the idea to use compression
algorithms to study and measure experimentally the kind of chaos
in intermittent dynamical systems can be practically exploited. In
\cite{licatone}, using a particular compression algorithm that is
suitable for null entropy strings, we performed experiments on
some examples of intermittent and weakly chaotic  dynamical
systems and obtained results that are close to the theoretical
predictions. Moreover, the study of weak chaos by compression
algorithms gives rise to new questions in data compression (the
search for algorithms that are optimal compressing zero entropy
strings).

We end remarking that we focused our interest onto maps of the
interval. A more general approach to define a weak chaos index
that is suitable for maps on a general metric space $X$ is to use
open covers instead of partitions (\cite{gal1},\cite{brudno}, see
also the remarks at the end of section 4 in \cite{licatone} ). In
this paper we chose to simplify the question by using admissible
partitions, however all the results given here for the partitions
hold also for the open covers. In  \cite{gal1} there is some
example of results of this kind for PL maps.

\subsection*{\emph{Acknowledgements}}
The authors wish to thank Stefano Isola for many suggestions and
useful discussions.

\end{document}